\documentclass[12pt]{article}
\usepackage{amsmath}
\usepackage{amssymb}
\usepackage{vatola}

\def\q{\quad}
\def\qq{\qquad}
\def\qtq#1{\q\t{#1}\q}
\def\mod#1{\ (\text{\rm mod}\ #1)}
\def\t{\text}
\def\f{\frac}
\def\e{\equiv}
\def\b{\binom}
\def\qp#1{q_p(#1)}
\def\sls#1#2{(\f{#1}{#2})}
 \def\ls#1#2{\big(\f{#1}{#2}\big)}
\def\Ls#1#2{\Big(\f{#1}{#2}\Big)}

\def\sumnp{\sum_{n=0}^{p-1}}
\def\sumkp{\sum_{k=0}^{p-1}}
\let \pro=\proclaim
\let \endpro=\endproclaim

\begin{document}
 \centerline {\bf
New congruences involving Ap\'ery-like numbers}
\par\q\newline
\centerline{Zhi-Hong Sun}\newline \centerline{School of Mathematics
and Statistics}\centerline{Huaiyin Normal University}
\centerline{Huaian, Jiangsu 223300, P.R. China} \centerline{Email:
zhsun@hytc.edu.cn} \centerline{Homepage:
http://maths.hytc.edu.cn/szh1.htm}
 \abstract{In this paper we present many congruences for several
 Ap\'ery-like sequences.
 \par\q
 }
 \endabstract
\section*{1. Introduction}

\par \par For $s>1$ let $\zeta(s)=\sum_{n=1}^{\infty}\f 1{n^s}$.
In 1979, in order to prove that $\zeta(2)$ and
 $\zeta(3)$ are irrational,
 Ap\'ery [Ap] introduced the Ap\'ery numbers $\{A_n\}$ and $\{A'_n\}$
 given by
 $$ A_n= \sum_{k=0}^n\binom nk^2\binom{n+k}k^2\qtq{and}A'_n=\sum_{k=0}^n\b
 nk^2\b{n+k}k.$$
It is well known  (see [B2]) that
$$\align &(n+1)^3A_{n+1}=(2n+1)(17n(n+1)+5)A_n-n^3A_{n-1}\q (n\ge 1),
\\&(n+1)^2A'_{n+1}=(11n(n+1)+3)A'_n+n^2A'_{n-1}\q (n\ge 1).\endalign$$
\par
Let $\Bbb Z$ and $\Bbb Z^+$ be the set of integers and the set of
positive integers, respectively.  The Ap\'ery-like numbers $\{u_n\}$
of the first kind satisfy
$$u_0=1, \ u_1=\ b,\ (n+1)^3u_{n+1}
=(2n+1)(an(n+1)+b)u_n-cn^3u_{n-1}\q (n\ge 1),\tag 1.1$$ where
$a,b,c\in\Bbb Z$ and $c\not=0$. Let $[x]$ be the greatest integer
not exceeding $x$, and let
$$\align
&D_n=\sum_{k=0}^n\b nk^2\b{2k}k\b{2n-2k}{n-k},\q T_n=\sum_{k=0}^n\b
nk^2\b{2k}n^2, \\& b_n=\sum_{k=0}^{[n/3]}\b{2k}k\b{3k}k\b
n{3k}\b{n+k}k(-3)^{n-3k}.
\endalign$$
Then $\{A_n\}$, $\{D_n\}$, $\{b_n\}$ and $\{T_n\}$  are Ap\'ery-like
numbers of the first kind with
$(a,b,c)=(17,5,1),(10,4,64),(-7,-3,81),(12,4,16)$, respectively. The
numbers $\{D_n\}$ are called Domb numbers, and $\{b_n\}$ are called
Almkvist-Zudilin numbers. For $\{A_n\}$, $\{D_n\}$, $\{b_n\}$ and
$\{T_n\}$ see A005259, A002895, A125143 and A290575 in Sloane's
database ``The On-Line Encyclopedia of Integer Sequences". For the
congruences involving $T_n$ see the author's recent paper [S10].
\par In 2009 Zagier [Z] studied the Ap\'ery-like numbers $\{u_n\}$
of the second kind  given by
$$u_0=1,\ u_1=b\qtq{and}(n+1)^2u_{n+1}=(an(n+1)+b)u_n-cn^2u_{n-1}\ (n\ge 1),\tag 1.2$$
where $a,b,c\in\Bbb Z$ and $c\not=0$. Let
 $$\align &f_n=\sum_{k=0}^n\b nk^3=\sum_{k=0}^n\b nk^2\b{2k}n,\q
 \\&S_n=\sum_{k=0}^{[n/2]}\b{2k}k^2\b n{2k}4^{n-2k}=\sum_{k=0}^n\b nk\b{2k}k\b{2n-2k}{n-k},
 \\&a_n=\sum_{k=0}^n\b nk^2\b{2k}k,
 \q Q_n=\sum_{k=0}^n\b nk(-8)^{n-k}f_k,
 \\&W_n=\sum_{k=0}^{[n/3]}\b{2k}k\b{3k}k\b n{3k}(-3)^{n-3k}.\endalign$$
   Zagier [Z] stated that $\{A'_n\},\ \{f_n\},\ \{S_n\},\ \{a_n\},\
 \{Q_n\}$
and $\{W_n\}$ are  Ap\'ery-like sequences of the second kind with
$(a,b,c)=(11,3,-1),(7,2,-8),(12,4,32),(10,3,9),(-17,$ $-6,72), (-9,$
$-3,27)$, respectively. The  sequence
 $\{f_n\}$ is called Franel numbers.
 In [JS,S8,S9,S11]  the author
 systematically investigated identities and congruences for sums
 involving $S_n$, $f_n$ or $W_n$.
  For $\{A'_n\},\ \{f_n\},\ \{S_n\},\ \{a_n\},\
 \{Q_n\}$
and $\{W_n\}$ see A005258, A000172, A081085, A002893, A093388 and
A291898 in Sloane's database ``The On-Line Encyclopedia of Integer
Sequences".
\par Ap\'ery-like numbers have fascinating properties and they are
concerned with Picard-Fuchs differential equation, modular forms,
hypergeometric series, elliptic curves, series for $\f 1{\pi}$,
supercongruences, binary quadratic forms, combinatorial identities,
Bernoulli numbers and Euler numbers. See for example [AO], [AZ],
[B1]-[B2], [CCS], [CZ], [CTYZ], [GMY], [SB] and [Z].

\par In this paper we mainly investigate identities and
congruences for two Ap\'ery-like sequences $\{G_n\}$ and $\{V_n\}$
originated from [AZ]. See A143583 and A036917 in Sloane's database
``The On-Line Encyclopedia of Integer Sequences". For
$n=0,1,2,\ldots$ let
$$G_n=\sum_{k=0}^n\b{2k}k^2\b{2n-2k}{n-k}4^{n-k},
\q V_n=\sum_{k=0}^n\b{2k}k^2\b{2n-2k}{n-k}^2.\tag 1.3$$ Using Maple
one can check that
$$\align &(n+1)^2G_{n+1}=(32n(n+1)+12)G_n-256n^2G_{n-1}\q (n\ge 1),\tag 1.4
\\&(n+1)^3V_{n+1}=(2n+1)(16n(n+1)+8)V_n-256n^3V_{n-1}\q(n\ge 1).\tag
1.5\endalign$$ Thus, $\{G_n\}$ is an Ap\'ery-like sequence of the
second kind with $(a,b,c)=(32,12,256)$, and $\{V_n\}$ is an
Ap\'ery-like sequence of the first kind with $(a,b,c)=(16,8,256)$.
 The first few values of $G_n$ and $V_n$ are shown
below:
$$\align &G_0=1,\ G_1=12,\ G_2=164,\ G_3=2352,\ G_4=34596,\ G_5=516912,
\ G_6=7806224,
\\&V_0=1,\ V_1=8,\ V_2=88,\ V_3=1088,\ V_4=14296,\ V_5=195008,
\ V_6=2728384.\endalign$$
 \par For
positive integers $a,b$ and $n$, if $n=ax^2+by^2$ for some integers
$x$ and $y$, we briefly write that $n=ax^2+by^2$. In Section 2 we
obtain new expressions and
 congruences for $G_n$. In particular, for any prime $p>3$,
$$\align &G_{p-1}\e (-1)^{\f{p-1}2}256^{p-1}\mod {p^2},\q
\sumnp\f{G_n}{16^n}\e 0\mod {p^2},
\\&G_{\f{p-1}2}\e \cases 4^px^2-2p\mod {p^2}& \t{if $p=x^2+4y^2\e
1\mod 4$,}
\\0\mod {p^2}&\t{if $p\e 3\mod 4$,}
\endcases
\\&\sum_{n=0}^{p-1}\f{\b{2n}nG_n}{128^n}
\e\cases 4x^2-2p\mod {p^2}&\t{if $p\e 1,3\mod 8$ and so
$p=x^2+2y^2$,}
\\0\mod {p^2}&\t{if $p\e 5,7\mod 8$,}
\endcases
\\&\sum_{n=0}^{p-1}\f{\b{2n}nnG_n}{128^n}\e 0\mod {p^2},
\\&\sum_{n=0}^{p-1}\f{\b{2n}nG_n}{48^n}
\e \cases 4x^2-2p\mod {p^2}&\t{if $p=x^2+3y^2\e 1\mod 3$,}
\\0\mod {p^2}&\t{if $p\e 2\mod 3$.}
\endcases
\endalign$$
 Section 3 is devoted to the properties of $V_n$. We show that
$$V_n =\sum_{k=0}^n\b nk\b{n+k}k(-16)^{n-k}G_k,$$
 and  for any prime $p>3$,
 $$\align &V_{p-1}\e 256^{p-1}\mod
 {p^3},\q\sum_{n=1}^{p-1}\f{V_n}{16^n}\e  \f 72p^3B_{p-3} \mod {p^4},
 \\&\sum_{n=0}^{p-1}\f{V_n}{(-16)^n}\e\cases 4x^2-2p \mod {p^2}&\t{if
$p\e 1\mod 4$ and so $p=x^2+4y^2$,} \\
0\mod p&\t{if $p\e 3\mod 4$,}\endcases
\endalign$$
where $\{B_n\}$ are Bernoulli numbers given by
$$B_0=1\qtq{and}\sum_{k=0}^{n-1}\b nkB_k=0\q(n\ge 2).$$
\par In Section 4 we prove several congruences involving
$b_n,A'_n,S_n$ and $Q_n$. In particular, we show that for any prime
$p>3$,
$$\align &b_{p-1}\e 81^{p-1}\mod{p^3},\q
\sum_{n=1}^{p-1}\f{nS_n}{8^n}\e 0\mod {p^2},
\\&\sum_{n=0}^{p-1}\f{Q_n}{(-8)^n}\e 1\mod {p^2},\q
\sum_{n=0}^{p-1}\f{Q_n}{(-9)^n}\e \Ls p3\mod {p^2},\endalign$$ where
$\ls ap$ is the Legendre symbol.

\par  For a prime $p$ let $\Bbb Z_p$ be the set of
rational numbers whose denominator is not divisible by $p$. For an
odd prime and $a\in\Bbb Z_p$ set $\qp a=(a^{p-1}-1)/p$. For $n\in
\Bbb Z^+$ define $H_n=1+\f 12+\cdots+\f 1n$ and $H_n^{(r)}=1+\f
1{2^r}+\cdots+\f 1{n^r}$ for $r\in\Bbb Z^+$. The Euler numbers
$\{E_n\}$ are defined by
$$ E_{2n-1}=0,\q E_0=1,\q E_{2n}=-\sum_{k=1}^{n}
\b {2n}{2k}E_{2n-2k}\q(n\ge 1).$$ For congruences involving $B_n$
and $E_n$ see [S1-S2].

\section*{2. Identities and congruences involving $G_n$}
\par Recall that
$$G_n=\sum_{k=0}^n\b{2k}k^2\b{2n-2k}{n-k}4^{n-k}\q(n=0,1,2,\ldots).$$

\pro{Theorem 2.1} For $n=0,1,2,\ldots$,
$$G_n=\sum_{k=0}^n\b nk(-1)^k\b{2k}k^216^{n-k}
.$$
\endpro
Proof. Set
$$G'_n=\sum_{k=0}^n\b nk(-1)^k\b{2k}k^216^{n-k}.$$
Then $G'_0=1=G_0$ and $G'_1=12=G_1$. Using sumtools in Maple we find
that
$$(n+1)^2G'_{n+1}=(32n(n+1)+12)G'_n-256n^2G'_{n-1}\q (n\ge 1).$$
Thus, applying (1.4) yields $G'_n=G_n$ for $n=0,1,2,\ldots$. This
proves the theorem.

\pro{Theorem 2.2} Let $p$ be a prime with $p>3$. Then
$$G_{p-1}\e (-1)^{\f{p-1}2}256^{p-1}+p^2\Big(E_{p-3}-8(-1)^{\f{p-1}2}\qp
2^2+\f 12\sum_{k=0}^{\f{p-1}2}\f{\b{2k}k^2}{16^k}H_k^2\Big) \mod
{p^3}.$$ and
$$G_{\f{p-1}2}\e \cases 4^px^2-2p\mod {p^2}&
\t{if $p=x^2+4y^2\e 1\mod 4$,}
\\0\mod {p^2}&\t{if $p\e 3\mod 4$.}
\endcases$$
\endpro
 Proof. By [S2, Lemma 2.9], for $1\le k\le p-1$,
 $$\b{p-1}k\e 1-pH_k+\f{p^2}2\big(H_k^2-H_k^{(2)}\big)\mod {p^3}.$$
 Thus, from Theorem 2.1 we derive that
$$G_{p-1}=\sum_{k=0}^{p-1}\b{p-1}k(-1)^k\b{2k}k^216^{p-1-k}
\e
16^{p-1}\sum_{k=0}^{p-1}\Big(1-pH_k+\f{p^2}2\big(H_k^2-H_k^{(2)}\big)\Big)
\f{\b{2k}k^2}{16^k}\mod {p^3}.$$ By [Su1, (1.7) and (1.9)],
$$\sum_{k=0}^{p-1}\f{\b{2k}k^2}{16^k}\e (-1)^{\f{p-1}2}-p^2E_{p-3}
\mod {p^3}.$$ By [Su7, Theorem 4.1] and [S1],
$$\align &\sum_{k=0}^{\f{p-1}2}\f{\b{2k}k^2}{16^k}H_k\e
2(-1)^{\f{p-1}2}H_{\f{p-1}2}\e (-1)^{\f{p-1}2}(-4\qp 2+2p\qp
2^2)\mod {p^2},
\\&\sum_{k=0}^{\f{p-1}2}\f{\b{2k}k^2}{16^k}H_k^{(2)}\e -4E_{p-3}\mod
p.\endalign$$ Thus,
$$\align G_{p-1}&\e
16^{p-1}\Big((-1)^{\f{p-1}2}-p^2E_{p-3}-(-1)^{\f{p-1}2}p(-4\qp
2+2p\qp 2^2)\\&\qq-\f
{p^2}2(-4E_{p-3})+\f{p^2}2\sum_{k=0}^{p-1}\f{\b{2k}k^2}{16^k}H_k^2\Big)
\\&\e 16^{p-1}\Big((-1)^{\f{p-1}2}(1+4p\qp 2)+p^2(E_{p-3}-2(-1)^{\f{p-1}2}\qp
2^2)+\f{p^2}2\sum_{k=0}^{p-1}\f{\b{2k}k^2}{16^k}H_k^2\Big)
\\&\e 16^{p-1}\Big((-1)^{\f{p-1}2}(1+p\qp
2)^4+p^2(E_{p-3}-8(-1)^{\f{p-1}2}\qp
2^2)+\f{p^2}2\sum_{k=0}^{p-1}\f{\b{2k}k^2}{16^k}H_k^2\Big)
\\&\e (-1)^{\f{p-1}2}256^{p-1}+p^2\Big(E_{p-3}-8(-1)^{\f{p-1}2}\qp
2^2+\f 12\sum_{k=0}^{\f{p-1}2}\f{\b{2k}k^2}{16^k}H_k^2\Big) \mod
{p^3}.
\endalign$$
 On the other hand, from Theorem 2.1 and [S3, Lemma 2.4],
$$\align G_{\f{p-1}2}&=\sum_{k=0}^{(p-1)/2}\b{\f{p-1}2}k(-1)^k\b{2k}k^216^{\f{p-1}2-k}
\e \sum_{k=0}^{(p-1)/2}\f{\b{2k}k}{4^k}\Big(1-p\sum_{i=1}^k\f
1{2i-1}\Big)\b{2k}k^216^{\f{p-1}2-k}
\\&=4^{p-1}\sum_{k=0}^{(p-1)/2}\f{\b{2k}k^3}{64^k}
\Big(1-p\sum_{i=1}^k\f 1{2i-1}\Big)\e
4^{p-1}\sum_{k=0}^{p-1}\f{\b{2k}k^3}{64^k}\Big(1-p\sum_{i=1}^k\f
1{2i-1}\Big)\mod{p^2}.
\endalign$$
It is well known (see [A]) that
$$\sum_{k=0}^{p-1}\f{\b{2k}k^3}{64^k}\e\cases 4x^2-2p\mod {p^2}
&\t{if $p\e 1\mod 4$ and so $p=x^2+4y^2$,}
\\0\mod {p^2}&\t{if $p\e 3\mod 4$.}
\endcases$$
By [T2, Theorem 5.2], $$\sum_{k=0}^{p-1}\f{\b{2k}k^3}{64^k}
\sum_{i=1}^k\f 1{2i-1}\e 0\mod p.$$ Thus,
$$G_{\f{p-1}2}\e 4^{p-1}\sum_{k=0}^{p-1}\f{\b{2k}k^3}{64^k}
\e\cases 4^{p-1}(4x^2-2p)\e 4^px^2-2p\mod {p^2}&\t{if $p=x^2+4y^2\e
1\mod 4$,}
\\0\mod {p^2}&\t{if $p\e 3\mod 4$.}
\endcases$$
This completes the proof.

\pro{Lemma 2.1 ([Su6, Lemma 2.2])} For $n=0,1,2,\ldots$,
$$\b{2n}nG_n=\sum_{k=0}^n\b{2k}k^2\b{4k}{2k}\b k{n-k}(-64)^{n-k}.$$
\endpro
\pro{Theorem 2.3} Let $p$ be an odd prime, $m\in\Bbb Z_p$ and
$m\not\e 0\mod p$. Then
$$\sum_{n=0}^{p-1}\f{\b{2n}nG_n}{m^n}
\e \sum_{k=0}^{(p-1)/2}\b{2k}k^2\b{4k}{2k}\Ls{m-64}{m^2}^k\e
\Big(\sum_{k=0}^{p-1}\f{\b{2k}k\b{4k}{2k}}{m^k}\Big)^2\mod {p^2}.$$
\endpro
Proof. Note that $p\mid \b{2k}k$ for $\f p2<k<p$. By Lemma 2.1, for
$c_0,c_1,\ldots,c_{p-1}\in\Bbb Z_p$,
$$\align &\sum_{n=0}^{p-1}\f{\b{2n}nc_nG_n}{m^n}\\&
=\sum_{n=0}^{p-1}c_n\sum_{k=0}^n\f{\b{2k}k^2\b{4k}{2k}}{(-64)^k}\b
k{n-k}\Ls{-64}m^n
=\sum_{k=0}^{p-1}\f{\b{2k}k^2\b{4k}{2k}}{(-64)^k}\sum_{n=k}^{p-1} \b
k{n-k}c_n\Ls{-64}m^n
\\&=\sum_{k=0}^{p-1}\f{\b{2k}k^2\b{4k}{2k}}{(-64)^k}\sum_{r=0}^{p-1-k}
\b krc_{k+r}\Ls{-64}m^{k+r}
=\sum_{k=0}^{p-1}\f{\b{2k}k^2\b{4k}{2k}}{m^k}\sum_{r=0}^{p-1-k} \b
krc_{k+r}\Ls{-64}m^{r}
\\&\e \sum_{k=0}^{(p-1)/2}\f{\b{2k}k^2\b{4k}{2k}}{m^k}\sum_{r=0}^k
\b krc_{k+r}\Ls{-64}m^{r}\mod {p^2}.\endalign$$ Thus, applying [S4,
Theorem 4.1],
$$\sum_{n=0}^{p-1}\f{\b{2n}nG_n}{m^n}
\e\sum_{k=0}^{(p-1)/2}\f{\b{2k}k^2\b{4k}{2k}}{m^k}\Big(1-\f
{64}m\Big)^k\e
\Big(\sum_{k=0}^{p-1}\f{\b{2k}k\b{4k}{2k}}{m^k}\Big)^2 \mod {p^2}.$$
This completes the proof.

\pro{Theorem 2.4} Let $p$ be an odd prime. Then
$$\align &\sum_{n=0}^{p-1}\f{\b{2n}nG_n}{128^n}
\e\cases 4x^2-2p\mod {p^2}&\t{if $p\e 1,3\mod 8$ and so
$p=x^2+2y^2$,}
\\0\mod {p^2}&\t{if $p\e 5,7\mod 8$,}
\endcases
\\&\sum_{n=0}^{p-1}\f{\b{2n}nnG_n}{128^n}\e 0\mod {p^2}.
\endalign$$
\endpro
Proof. Taking $m=128$ in Theorem 2.3 and then applying [S4, Theorem
4.3, (1.6) and Lemma 2.2] yields the first congruence. By the proof
of Theorem 2.3,
$$ \sum_{n=0}^{p-1}\f{\b{2n}nnG_n}{128^n}
\e\sum_{k=0}^{(p-1)/2}\f{\b{2k}k^2\b{4k}{2k}}{128^k}\sum_{r=0}^k \b
kr \f {k+r}{(-2)^r} \mod {p^2}.$$ Note that
$$\sum_{r=0}^k\b kr\f{k+r}{(-2)^r}=k\sum_{r=0}^k\b kr\f 1{(-2)^r}-\f
k2\sum_{r=1}^k\b{k-1}{r-1}\f 1{(-2)^{r-1}}=\f k{2^k}-\f k2\cdot \f
1{2^{k-1}}=0.$$ We then obtain the second congruence.

\pro{Theorem 2.5} Let $p>3$ be a prime. Then
$$\align &\sum_{n=0}^{p-1}\f{\b{2n}nG_n}{72^n}
\e \cases 4x^2-2p\mod {p^2}&\t{if $p=x^2+4y^2\e 1\mod 4$,}
\\0\mod {p^2}&\t{if $p\e 3\mod 4$,}
\endcases
\\&\sum_{n=0}^{p-1}\f{\b{2n}nG_n}{576^n} \e\cases 4x^2\mod p&\t{if
$p=x^2+4y^2\e 1\mod 4$,}
\\0\mod {p^2}&\t{if $p\e 3\mod 4$,}
\endcases
\\&\sum_{n=0}^{p-1}\f{\b{2n}nG_n}{48^n}
\e \cases 4x^2-2p\mod {p^2}&\t{if $p=x^2+3y^2\e 1\mod 3$,}
\\0\mod {p^2}&\t{if $p\e 2\mod 3$,}
\endcases
\\&\sum_{n=0}^{p-1}\f{\b{2n}nG_n}{(-192)^n} \e\cases 4x^2\mod
p&\t{if $p=x^2+3y^2\e 1\mod 3$,}
\\0\mod {p^2}&\t{if $p\e 2\mod 3$,}
\endcases
\\&\sum_{n=0}^{p-1}\f{\b{2n}nG_n}{63^n}
\e \sum_{n=0}^{p-1}\f{\b{2n}nG_n}{(-4032)^n} \e\cases 4x^2\mod
p&\t{if $p=x^2+7y^2\e 1,2,4\mod 7$,}
\\0\mod {p^2}&\t{if $p\e 3,5,6\mod 7$.}
\endcases
\endalign$$
\endpro
Proof. This is immediate from Theorem 2.3, [S4, Theorems 2.2-2.4]
and [WS, Remarks 5.1 and 5.2].

\pro{Theorem 2.6} Let $p$ be an odd prime. Then
$$\sum_{n=0}^{p-1}\f{G_n}{16^n}\e p^2\Big(1-
\sum_{k=0}^{\f{p-1}2}\f{\b{2k}k^2}{16^k(k+1)}H_k\Big)\mod {p^3}.$$
\endpro
Proof. From [G, (1.52)],
$$\sum_{n=r}^m\b nr=\b{m+1}{r+1}.\tag 2.1$$
By Theorem 2.1 and (2.1),
$$\align \sum_{n=0}^{p-1}\f{G_n}{16^n}
&=\sum_{n=0}^{p-1}\sum_{k=0}^{p-1}\b nk(-1)^k\f{\b{2k}k^2}{16^k}
=\sum_{k=0}^{p-1}\f{\b{2k}k^2}{(-16)^k}\sum_{n=k}^{p-1}\b nk
\\&=\sum_{k=0}^{p-1}\f{\b{2k}k^2}{(-16)^k}\b p{k+1}
=\sum_{k=0}^{p-1}\f{\b{2k}k^2}{(-16)^k}\cdot\f p{k+1}\b{p-1}k
\\&\e \sum_{k=0}^{\f{p-1}2}\f{\b{2k}k^2}{(-16)^k}\cdot\f p{k+1}\b{p-1}k
+\f{\b{2(p-1)}{p-1}^2}{16^{p-1}}
\\&\e p\sum_{k=0}^{\f{p-1}2}\f{\b{2k}k^2}{16^k(k+1)}(1-pH_k)+\f
{p^2}{(2p-1)^2}\cdot\f{\b{2p-1}{p-1}^2}{16^{p-1}}
\\&\e p\sum_{k=0}^{\f{p-1}2}\f{\b{2k}k^2}{16^k(k+1)}
-p^2\sum_{k=0}^{\f{p-1}2}\f{\b{2k}k^2}{16^k(k+1)}H_k+p^2
 \mod {p^3}.\endalign$$
 From [G, (3.100)] we know that
$$\sum_{k=0}^n\b{2k}k(-1)^k\b{n+k}{2k}\f {1}{k+b}
=(-1)^n\f{(b-1)(b-2)\cdots(b-n)}{b(b+1)\cdots(b+n)}.\tag 2.2$$
 From [S3, Lemma 2.2],
$$\f{\b{2k}k}{(-16)^k}\e \b{\f{p-1}2+k}{2k}
\mod{p^2}\qtq{for}k=1,2,\ldots,\f{p-1}2.\tag 2.3$$ Appealing to
(2.2) and (2.3),
$$\sum_{k=0}^{\f{p-1}2}\f{\b{2k}k^2}{16^k}\cdot\f p{k+1}
\e \sum_{k=0}^{\f{p-1}2}\b{2k}k(-1)^k\b{\f{p-1}2+k}{2k} \f p{k+1}=0
\mod{p^3}.$$ Combining the above prove the theorem.

\pro{Theorem 2.7} Let $p$ be an odd prime, $m\in\Bbb Z_p$, $m\not\e
0\mod p$ and $t=1-\f{32}m$. Then
$$\sum_{n=0}^{p-1}\f{G_n}{m^n}
\e\Ls{m(m-16)}p\sum_{k=0}^{\f{p-1}2}\f{\b{2k}k^2}{m^k}\e
-\Ls{-3(1+t)}p\sum_{x=0}^{p-1}\Ls{x^3-3(t^2+3)x+2t(t^2-9)}p\mod p.$$
\endpro
Proof. Clearly,
$$\align \sum_{n=0}^{p-1}\f{G_n}{m^n}&
=\sum_{n=0}^{p-1}\sum_{k=0}^n\b{2k}k^2\b{2(n-k)}{n-k}\f{4^{n-k}}{m^n}
=\sum_{k=0}^{p-1}\f{\b{2k}k^2}{m^k}\sum_{n=k}^{p-1}\b{2(n-k)}{n-k}
\Ls 4m^{n-k}
\\&=\sum_{k=0}^{p-1}\f{\b{2k}k^2}{m^k}\sum_{r=0}^{p-1-k}
\b{2r}r\Ls
4m^r=\sum_{k=0}^{p-1}\f{\b{2k}k^2}{m^k}\sum_{r=0}^{p-1-k}\b{-\f 12}r
\Big(-\f{16}m\Big)^r
\\&\e \sum_{k=0}^{(p-1)/2}\f{\b{2k}k^2}{m^k}\sum_{r=0}^{p-1-k}\b{-\f
12}r \Big(-\f{16}m\Big)^r \e
\sum_{k=0}^{(p-1)/2}\f{\b{2k}k^2}{m^k}\sum_{r=0}^{(p-1)/2}\b{\f
{p-1}2}r \Big(-\f{16}m\Big)^r
\\&=\Big(1-\f{16}m\Big)^{\f{p-1}2}\sum_{k=0}^{(p-1)/2}\f{\b{2k}k^2}{m^k}
\e \Ls{m(m-16)}p\sum_{k=0}^{(p-1)/2}\f{\b{2k}k^2}{m^k}\mod p.
\endalign$$
Now applying [S3, Theorem 2.11] yields the result.
 \pro{Theorem
2.8} Let $p$ be an odd prime. Then
$$\align &\sum_{n=0}^{p-1}\f{G_n}{(-16)^n}\e (-1)^{\f{p-1}4}
\sum_{n=0}^{p-1}\f{G_n}{8^n}\e
(-1)^{\f{p-1}4}\sum_{n=0}^{p-1}\f{G_n}{32^n} \\& \e\cases 2x\mod
p&\t{if
$4\mid p-1$ and so $p=x^2+4y^2$ with $4\mid x-1$,} \\
0\mod p&\t{if $p\e 3\mod 4$.}\endcases\endalign$$
\endpro
Proof. Taking $m=-16,8,32$ in Theorem 2.7 and then applying [S3,
Theorems 2.2,2.9 and Corollary 2.3] yields the result.
\par\q

\pro{Conjecture 2.1} For any prime $p>3$,
$$G_{p-1}\e (-1)^{\f{p-1}2}256^{p-1}+3p^2E_{p-3}\mod{p^3}.$$
\endpro
\pro{Conjecture 2.2} For any odd prime $p$,
$$\sum_{n=0}^{p-1}\f{G_n}{16^n}\e \big(4(-1)^{\f{p-1}2}-3\big)p^2
\mod {p^3}.$$
\endpro

\pro{Conjecture 2.3} Let $p$ be an odd prime. Then
$$\align &\sum_{n=0}^{p-1}\f{G_n}{(-16)^n}\e\cases 2x-\f p{2x}\mod
{p^2}&\t{if $4\mid p-1$, $p=x^2+4y^2$ and $4\mid x-1$,}\\
p\b{\f{p-3}2}{\f{p-3}4}^{-1}\mod {p^2}&\t{if $p\e 3\mod
4$,}\endcases
\\&\sum_{n=0}^{p-1}\f{G_n}{8^n}\e
\cases (-1)^{\f{p-1}4}(2x-\f p{2x})\mod
{p^2}&\t{if $4\mid p-1$, $p=x^2+4y^2$ and $4\mid x-1$,}\\
 (-1)^{\f{p+1}4}\f 32
p\b{\f{p-3}2}{\f{p-3}4}^{-1}\mod {p^2}&\t{if $p\e 3\mod
4$,}\endcases
\\&\sum_{n=0}^{p-1}\f{G_n}{32^n}\e
\cases (-1)^{\f{p-1}4}(2x-\f p{2x})\mod
{p^2}&\t{if $4\mid p-1$, $p=x^2+4y^2$ and $4\mid x-1$,}\\
 (-1)^{\f{p-3}4}\f
p2\b{\f{p-3}2}{\f{p-3}4}^{-1}\mod {p^2}&\t{if $p\e 3\mod
4$.}\endcases
\endalign$$
\endpro

\pro{Conjecture 2.4} Let $p>3$ be a prime. Then
$$\align &\sum_{n=0}^{p-1}\f{\b{2n}nnG_n}{72^n}
\e -8\sum_{n=0}^{p-1}\f{\b{2n}nnG_n}{576^n} \e\cases 4x^2-3p\mod
{p^2}\\\qq\qq\qq\t{if $p=x^2+4y^2\e 1\mod 4$,}
\\p\mod {p^2}\q\ \t{if $p\e 3\mod 4$,}
\endcases
\\&\sum_{n=0}^{p-1}\f{\b{2n}nnG_n}{48^n}
\e 4\sum_{n=0}^{p-1}\f{\b{2n}nnG_n}{(-192)^n}\e\cases 3p-4x^2\mod
{p^2}\\\qq\qq\qq\q\t{if $p=x^2+3y^2\e 1\mod 3$,}
\\-p\mod {p^2}\q\ \t{if $p\e 2\mod 3$,}
\endcases
\\&\sum_{n=0}^{p-1}\f{\b{2n}nnG_n}{63^n}
\e 64\sum_{n=0}^{p-1}\f{\b{2n}nnG_n}{(-4032)^n} \e\cases
-8(4x^2-3p)\mod {p^2}\\\qq\qq\t{if $p=x^2+7y^2\e 1,2,4\mod 7$,}
\\-8p\mod {p^2}\q\ \t{if $p\e 3,5,6\mod 7$}
\endcases
\endalign$$
and
$$\sum_{n=0}^{p-1}\f{\b{2n}nG_n}{128^n}
\e \cases 4x^2-2p-\f{p^2}{4x^2}\mod {p^3}& \t{if $p=x^2+2y^2 \e
1,3\mod 8$,}\\\f 13p^2\b{[p/4]}{[p/8]}^{-2}\mod {p^3}&\t{if $p\e
5\mod 8$,}
\\-\f 32p^2\b{[p/4]}{[p/8]}^{-2}\mod {p^3}&\t{if $p\e
7\mod 8$.}\endcases$$
\endpro

\pro{Conjecture 2.5} Let $p$ be an odd prime and $m,r\in\Bbb Z^+$.
Then
$$G_{mp^r}\e G_{mp^{r-1}}\mod {p^{2r}}.$$
\endpro

\pro{Conjecture 2.6} Suppose that $p$ is a prime with $p\e 3\mod 4$
, $m\in\{1,3,5,\ldots\}$ and $r\in\{2,3,4,\ldots\}$. Then
$$G_{\f{mp^r-1}2}\e p^2G_{\f{mp^{r-2}-1}2}\mod {p^{2r-1}}.$$
\endpro
\pro{Conjecture 2.7} Suppose that $p$ is a prime of the form $4k+1$,
 $p=x^2+4y^2$ , $m\in\{1,3,5,\ldots\}$ and $r\in\{2,3,4,\ldots\}$.
Then
$$G_{\f{mp^r-1}2}\e (4x^2-2p)G_{\f{mp^{r-1}-1}2}-
p^2G_{\f{mp^{r-2}-1}2}\mod {p^r}.$$

\section*{3. Congruences for $\{V_n\}$}
Recall that
$$V_n=\sum_{k=0}^n\b{2k}k^2\b{2n-2k}{n-k}^2.$$
By [W],
$$V_n=\sum_{k=0}^n\b{2k}k\b nk\b{n+k}{k}(-1)^k16^{n-k}=
\sum_{k=0}^n\b{2k}k^3\b{n+k}{2k}(-1)^k16^{n-k}.\tag 3.1$$

 \pro{Theorem 3.1} For any nonnegative integer $n$,
 $$V_n=\sum_{k=0}^n\b nk\b{n+k}k(-16)^{n-k}G_k.$$
 \endpro
 Proof. By Theorem 2.1,
$G_n=\sum_{k=0}^n\b nk(-1)^k\b{2k}k^216^{n-k}$. Applying the
binomial inversion formula,
$$\sum_{k=0}^n\b nk(-1)^k\f{G_k}{16^k}=\f{\b{2n}n^2}{16^n}.$$
Hence, from [S7, Theorem 2.2] we deduce
$$\align\sum_{k=0}^n\b nk\b{n+k}k\f{G_k}{(-16)^k}
&=\sum_{k=0}^n\b nk\b{n+k}k(-1)^{n-k}\sum_{s=0}^k\b
ks\f{G_s}{(-16)^s} \\&=\sum_{k=0}^n\b
nk\b{n+k}k(-1)^{n-k}\f{\b{2k}k^2}{16^k}
\\&=\sum_{k=0}^n\b{2k}k^3\b{n+k}{2k}\f{(-1)^n}{(-16)^k}
.\endalign$$ Thus, appealing to (3.1),
$$\sum_{k=0}^n\b nk\b{n+k}k(-16)^{n-k}G_k
=\sum_{k=0}^n\b{2k}k^3\b{n+k}{2k}(-1)^k16^{n-k}=V_n.$$ This
completes the proof.

\pro{Theorem 3.2} For $|x|<\f 1{16}$,
$$\sum_{n=0}^{\infty}V_nx^n=\f 1{1-16x}\sum_{k=0}^{\infty}\b{2k}k^3\Ls
{-x}{(1-16x)^2}^k.$$
\endpro
Proof. By (3.1), for $|x|<\f 1{16}$,
$$\align&\sum_{n=0}^{\infty}V_nx^n
\\&=\sum_{n=0}^{\infty}\sum_{k=0}^n\b{2k}k^3\b{n+k}{2k}(-1)^k16^{n-k}x^n
=\sum_{k=0}^{\infty}\b{2k}k^3(-x)^k\sum_{n=k}^{\infty}\b{n+k}{n-k}(16x)^{n-k}
\\&=\sum_{k=0}^{\infty}\b{2k}k^3(-x)^k\sum_{r=0}^{\infty}\b{2k+r}r(16x)^r
=\sum_{k=0}^{\infty}\b{2k}k^3(-x)^k\sum_{r=0}^{\infty}\b{-2k-1}r(-16x)^r
\\&=\sum_{k=0}^{\infty}\b{2k}k^3(-x)^k(1-16x)^{-2k-1}
=\f 1{1-16x}\sum_{k=0}^{\infty}\b{2k}k^3\Ls
{-x}{(1-16x)^2}^k.\endalign$$

\pro{Theorem 3.3} Let $p$ be a prime with $p>3$. Then
$$V_{p-1}\e 256^{p-1}\mod {p^3}.$$
\endpro
Proof. Note that $p\mid \b{2k}k$ for $\f p2<k<p$. From (3.1),
$$\align V_{p-1}&=\sum_{k=0}^{p-1}\b{2k}k^3\b{p-1+k}{2k}(-1)^k16^{p-1-k}
\e 16^{p-1}\sum_{k=0}^{(p-1)/2}\f{\b{2k}k^3}{(-16)^k}\b{p-1+k}{2k}
\\&= 16^{p-1}+16^{p-1}\sum_{k=1}^{(p-1)/2}\f{\b{2k}k^3}{(-16)^k}\cdot\f
{p(p-k)(p^2-1^2)(p^2-2^2)\cdots(p^2-(k-1)^2)}{(2k)!}
 \\&\e  16^{p-1}+16^{p-1}\sum_{k=1}^{(p-1)/2}\f{\b{2k}k^3}{(-16)^k}
 \cdot\f{p(p-k)(-1)^{k-1}k!^2}{k^2\cdot(2k)!}
 \\&=16^{p-1}\Big(1-p^2\sum_{k=1}^{(p-1)/2}\f{\b{2k}k^2}{16^kk^2}
 +p\sum_{k=1}^{(p-1)/2}\f{\b{2k}k^2}{16^kk}\Big)\mod {p^3}.
 \endalign$$
 By [T1] or [S10, p.1692],
 $$\sum_{k=1}^{(p-1)/2}\f{\b{2k}k^2}{16^kk^2}\e -8q_p(2)^2\pmod p,\q
\sum_{k=1}^{(p-1)/2}\f{\b{2k}k^2}{16^kk}\e
4q_p(2)-2pq_p(2)^2\pmod{p^2}.$$ Hence
$$\align V_{p-1}&\e 16^{p-1}\big(1-p^2(-8\qp
2^2)+p(4q_p(2)-2pq_p(2)^2)\big)=16^{p-1}\big(1+4p\qp 2+6p^2\qp
2^2\big) \\&\e 16^{p-1}(1+p\qp 2)^4=256^{p-1}\mod {p^3}\endalign$$
as claimed.

\pro{Theorem 3.4} Let $p$ be an odd prime, $m\in\Bbb Z_p$ and
$m\not\e 0,\pm 16\mod p$. Then
$$\sumnp\f{V_n}{m^n}\e \sumkp\b{2k}k\f{G_k}{(32+m+256/m)^k}
\e\sum_{k=0}^{\f{p-1}2}\f{\b{2k}k^3}{(32-m-256/m)^k}\mod p.$$
\endpro
Proof. By (3.1),
$$\align \sumnp\f{V_n}{m^n}&=\sumnp\sum_{k=0}^n\b{2k}k^3
\b{n+k}{2k}(-1)^k\f {16^{n-k}}{m^n} =\sum_{k=0}^{p-1}
\f{\b{2k}k^3}{(-m)^k}\sum_{n=k}^{p-1}\b{n+k}{n-k}\Ls {16}m^{n-k}
\\&=\sum_{k=0}^{p-1}
\f{\b{2k}k^3}{(-m)^k}\sum_{r=0}^{p-1-k}\b{2k+r}{r}\Ls{16}m^r
=\sum_{k=0}^{p-1}
\f{\b{2k}k^3}{(-m)^k}\sum_{r=0}^{p-1-k}\b{-2k-1}{r}\Ls{-16}m^r
\\&\e \sum_{k=0}^{(p-1)/2}
\f{\b{2k}k^3}{(-m)^k}\sum_{r=0}^{p-1-k}\b{p-1-2k}{r}\Ls{-16}m^r
=\sum_{k=0}^{(p-1)/2}
\f{\b{2k}k^3}{(-m)^k}\Big(1-\f{16}m\Big)^{p-1-2k}
\\&\e \sum_{k=0}^{(p-1)/2}\b{2k}k^3\Ls{1}{-m(1-16/m)^2}^k
=\sum_{k=0}^{(p-1)/2}\f{\b{2k}k^3}{(32-m-256/m)^k}\mod p.
\endalign$$
On the other hand, from [S8, Lemma 2.4], for $u\not\e 1\mod p$,
$$\sum_{n=0}^{p-1}u^n\sum_{k=0}^n \b nk\b{n+k}kc_k\e
\sum_{k=0}^{p-1}\b{2k}k\Ls u{(1-u)^2}^kc_k\mod p.$$ Since
$\f{V_n}{(-16)^n}=\sum_{k=0}^n\b nk\b{n+k}k\f{G_k}{(-16)^k}$ by
Theorem 3.1, taking $c_k=\f{G_k}{(-16)^k}$ and $u=-\f{16}m$ yields
$$\align \sum_{n=0}^{p-1}\f{V_n}{m^n}&=\sum_{n=0}^{p-1}\Ls{-16}m^n\sum_{k=0}^n
\b nk\b{n+k}k\f{G_k}{(-16)^k} \e \sum_{k=0}^{p-1}\b{2k}k
\Big(\f{-16/m}{(1+16/m)^2}\Big)^k\f{G_k}{(-16)^k}
\\&=\sum_{k=0}^{p-1}\b{2k}k\f {G_k}{(m+32+256/m)^k}\mod p.\endalign$$
This proves the theorem.

\pro{Theorem 3.5} Let $p$ be an odd prime. Then
$$\sum_{n=0}^{p-1}\f{V_n}{8^n}\e
\sum_{n=0}^{p-1}\f{V_n}{32^n}\e \cases 4x^2\mod p&\t{if
$p=x^2+4y^2\e 1\mod 4$,}
\\0\mod p&\t{if $p\e 3\mod 4$}
\endcases$$
and
$$\sum_{n=0}^{p-1}\f{V_n}{(-16)^n}\e\cases 4x^2-2p \mod {p^2}&\t{if
$p\e 1\mod 4$ and so $p=x^2+4y^2$,} \\
0\mod p&\t{if $p\e 3\mod 4$.}\endcases$$
\endpro
Proof. Taking $m=8$ and $32$ in Theorem 3.4 gives
$$\sum_{n=0}^{p-1}\f{V_n}{8^n}\e
\sum_{n=0}^{p-1}\f{V_n}{32^n}\e\sum_{k=0}^{p-1}\b{2k}k\f{G_k}{72^k}
\mod p.$$ Now applying Theorem 2.5 yields the first congruence.
 \par
By Theorem 3.1, $\f{V_n}{(-16)^n}=\sum_{k=0}^n\b
nk\b{n+k}k\f{G_k}{(-16)^k}$. Taking $c_k=\f{G_k}{(-16)^k}$ in [S9,
Lemma 2.2] gives
$$\align \sum_{n=0}^{p-1}\f{V_n}{(-16)^n}
&\e\sum_{k=0}^{p-1}\f p{2k+1}\cdot\f {G_k}{16^k}
\\&=\f{G_{\f{p-1}2}}{16^{\f{p-1}2}}+p\sum_{k=0}^{\f{p-3}2}\Big(\f{G_k}{16^k(2k+1)}
+\f{G_{p-1-k}}{16^{p-1-k}(2(p-1-k)+1)}\Big)\mod {p^3}.\endalign$$ By
[S11, Theorem 3.1] and Theorem 2.2,
$$G_n\e G_{p-1}256^nG_{p-1-n}\e
(-1)^{\f{p-1}2}256^nG_{p-1-n}\mod p\qtq{for}n=0,1,\ldots,p-1.$$
Hence
$$\align \sum_{n=0}^{p-1}\f{V_n}{(-16)^n}
&\e \f{G_{\f{p-1}2}}{16^{\f{p-1}2}}
+p\sum_{k=0}^{\f{p-3}2}\Big(\f{G_k}{16^k(2k+1)}
-\f{256^kG_{p-1-k}}{16^k(2k+1)}\Big)
\\&\e \f{G_{\f{p-1}2}}{4^{p-1}}
+p\sum_{k=0}^{\f{p-3}2}\Big(\f{G_k}{16^k(2k+1)}
-(-1)^{\f{p-1}2}\f{G_k}{16^k(2k+1)}\Big)\mod {p^2}.\endalign$$ Now
applying Theorem 2.2 gives
$$ \sum_{n=0}^{p-1}\f{V_n}{(-16)^n}
\e\cases \f{G_{\f{p-1}2}}{4^{p-1}}\e \f{4^px^2-2p}{4^{p-1}}\e
4x^2-2p\mod{p^2}&\t{if $p=x^2+4y^2\e 1\mod 4$,}
\\\f{G_{\f{p-1}2}}{4^{p-1}}\e 0\mod p&\t{if $p\e 3\mod 4$.}
\endcases$$
This completes the proof.

\pro{Theorem 3.6} Let $p$ be an odd prime. Then
$$\sum_{n=1}^{p-1}\f{V_n}{16^n}
\e \f 72p^3B_{p-3} \mod {p^4}.$$
\endpro
Proof. By (3.1),
$$\align\sum_{n=0}^{p-1}\f{V_n}{16^n}
&=\sum_{n=0}^{p-1}\sum_{k=0}^n\b{2k}k^3\b{n+k}{2k}\f 1{(-16)^k}
=\sum_{k=0}^{p-1}\f{\b{2k}k^3}{(-16)^k}\sum_{n=k}^{p-1}\b {n+k}{2k}
\\&=\sum_{k=0}^{p-1}\f{\b{2k}k^3}{(-16)^k}\b{p+k}{2k+1}
=\sum_{k=0}^{p-1}\f{\b{2k}k^3}{(-16)^k}\cdot\f
p{(2k+1)!}(p^2-1^2)\cdots(p^2-k^2) \\&\e
\sum_{k=0}^{p-1}\f{\b{2k}k^3}{(-16)^k}\cdot\f p{(2k+1)!}
(-1^2)(-2^2)\cdots(-k^2)(1-p^2H_k^{(2)})
\\&=\sum_{k=0}^{p-1}\f{\b{2k}k^2}{16^k}\cdot\f
p{2k+1}(1-p^2H_k^{(2)})\mod{p^4}.
\endalign$$
By [Su5, Theorem 1.2],
$$\align&\sum_{k=0}^{\f{p-3}2}\f{\b{2k}k^2}{(2k+1)16^k}
\e -2\qp 2-p\qp 2^2+\f 5{12}p^2B_{p-3}\mod {p^3},
\\&\sum_{k=\f{p+1}2}^{p-1}\f{\b{2k}k^2}{(2k+1)16^k}\e -\f
74p^2B_{p-3}\mod {p^3}.\endalign$$ By [S1, Theorem 5.2 and Corollary
5.2],
$$\align&H_{\f{p-1}2}\e -2\qp 2+p\qp 2^2-\f 23p^2\qp 2^3-\f
7{12}p^2B_{p-3}\mod {p^3},
\\&H_{\f{p-1}2}^{(2)}\e \f 73pB_{p-3}\mod {p^2},\q
H_{\f{p-1}2}^{(3)}\e -2B_{p-3}\mod p.\endalign$$ Appealing to [S2,
Lemma 2.9],
$$\align (-1)^{\f{p-1}2}\b{p-1}{\f{p-1}2}
&\e
1-pH_{\f{p-1}2}+\f{p^2}2\big(H_{\f{p-1}2}^2-H_{\f{p-1}2}^{(2)}\big)
-\f{p^3}6\big(H_{\f{p-1}2}^3-3H_{\f{p-1}2}H_{\f{p-1}2}^{(2)}+2H_{\f{p-1}2}^{(3)}\big)
\\&\e 1+2p\qp 2-p^2\qp 2^2+\f 23p^3\qp 2^3+\f 7{12}p^3B_{p-3}
\\&\q+\f {p^2}2\Big((-2\qp 2+p\qp 2^2)^2-\f 73pB_{p-3}\Big)
-\f{p^3}6\big((-2\qp 2)^3+2(-2B_{p-3})\big)
\\&\e 1+2p\qp 2+p^2\qp 2^2+\f{p^3}{12}B_{p-3}
=4^{p-1}+\f{p^3}{12}B_{p-3}\mod {p^4}.\endalign$$ Thus,
$$\f{\b{p-1}{\f{p-1}2}^2}{16^{\f{p-1}2}}\big(1-p^2H_{\f{p-1}2}^{(2)}\big)
\e\f {(4^{p-1}+\f{p^3}{12}B_{p-3})^2}{4^{p-1}}\Big(1-p^2\cdot\f
73pB_{p-3}\Big)\e 4^{p-1}-\f{13}6p^3B_{p-3}\mod {p^4}.$$ Now, from
the above we deduce that
$$\aligned \sum_{n=0}^{p-1}\f{V_n}{16^n}&
\e\sum\Sb k=0\\k\not=\f{p-1}2\endSb^{p-1} \f{\b{2k}k^2}{16^k}
\cdot\f p{2k+1}\big(1-p^2H_k^{(2)}\big) +
\f{\b{p-1}{\f{p-1}2}^2}{16^{\f{p-1}2}}\big(1-p^2H_{\f{p-1}2}^{(2)}\big)
\\&\e p\Big(-2\qp 2-p\qp 2^2+\f 5{12}p^2B_{p-3}-\f 74p^2B_{p-3}\Big)
\\&\qq-p^3\sum\Sb k=0\\k\not=\f{p-1}2\endSb^{p-1}
\f{\b{2k}k^2}{16^k(2k+1)}H_k^{(2)}+4^{p-1}-\f{13}6p^3B_{p-3}
\\&\e 1-p^3\Big(\sum_{k=0}^{\f{p-3}2}\f{\b{2k}k^2}{16^k(2k+1)}H_k^{(2)}+
\f 72B_{p-3}\Big) \mod {p^4}.\endaligned$$
 In 2019, Mao, Wang and Wang [MWW]
proved Z.W. Sun's conjecture:
$$\sum_{k=0}^{\f{p-3}2}\f{\b{2k}k^2}{16^k(2k+1)}H_k^{(2)}
\e -7B_{p-3}\mod p.$$ Thus, $$\sum_{n=0}^{p-1}\f{V_n}{16^n} \e 1-p^3
\Big(-7B_{p-3}+\f 72B_{p-3}\Big)\e 1+\f 72p^3B_{p-3}\mod {p^4}.$$
 This proves the theorem.
 \par\q
 \par{\bf Remark 3.1} Let $p$ be an odd prime.
In [Su6], Z.W. Sun conjectured that
$$\sum_{n=1}^{p-1}\f{nV_n}{32^n}\e -2p^3E_{p-3}\mod {p^4}\tag 3.2$$
and proved the congruence modulo $p^3$. (3.2) was solved by Mao and
Cao in [MC]. In [Su2], Z.W. Sun conjectured that
$$\align
&\sum_{n=0}^{p-1}(n+1)\f{V_n}{8^n}\e (-1)^{\f{p-1}2}p+5p^3
E_{p-3}\mod {p^4},\tag 3.3
\\&\sum_{n=0}^{p-1}(2n+1)\f{V_n}{(-16)^n}\e (-1)^{\f{p-1}2}p+3p^3
E_{p-3}\mod {p^4}.\tag 3.4\endalign$$ These congruences were
recently proved by Wang [W].
\par\q
\par Based on calculations by Maple, we pose the following
conjectures.
 \pro{Conjecture 3.1} Let $p$ be a prime with $p>3$.
Then
$$V_{p-1}\e 256^{p-1}-\f 32p^3B_{p-3}\mod {p^4}.$$
\endpro

 \pro{Conjecture
3.2} Let $p$ be a prime with $p>3$. Then
$$\sum_{n=0}^{p-1}\f{V_n}{8^n}
\e \sum_{n=0}^{p-1}\f{V_n}{(-16)^n} \e\cases
4x^2-2p-\f{p^2}{4x^2}\mod {p^3}&\t{if $p=x^2+4y^2\e 1\mod 4$,}
\\\f 34p^2\b{\f{p-3}2}{\f{p-3}4}^{-2}\mod {p^3}
&\t{if $p\e 3\mod 4$}\endcases$$ and
$$\sum_{n=0}^{p-1}\f{V_n}{32^n}
 \e\cases
4x^2-2p-\f{p^2}{4x^2}\mod {p^3}&\t{if $p=x^2+4y^2\e 1\mod 4$,}
\\-\f 14p^2\b{\f{p-3}2}{\f{p-3}4}^{-2}\mod {p^3}
&\t{if $p\e 3\mod 4$.}\endcases$$
\endpro

\pro{Conjecture 3.3} Let $p>3$ be a prime and $m,r\in\Bbb Z^+$. Then
$$V_{mp^r}\e V_{mp^{r-1}}+p^{3r}\cdot\f{V_{mp}-V_m}{p^3}\mod
{p^{3r+1}}.$$
\endpro

\pro{Conjecture 3.4} Let $p$ be an odd prime and $m,r\in\Bbb Z^+$.
Then
$$ V_{mp^r-1}\e 256^{mp^{r-1}(p-1)}V_{mp^{r-1}-1}\mod {p^{3r}}.$$
Moreover, for $p>3$,
$$ \f{V_{mp^r-1}-256^{mp^{r-1}(p-1)}V_{mp^{r-1}-1}}{p^{3r}}
\e \f{V_{mp-1}-256^{m(p-1)}V_{m-1}}{p^3}\mod p.$$
\endpro

\section*{4. Congruences for $b_n,A_n,S_n$ and $Q_n$}

\par
Let $\{b_n\}$ be the  Almkvist-Zudilin numbers. From
  [S6, (5.1) and (5.2)] or [CZ, Corollary 4.3] we know that
$$\aligned
b_n&=\sum_{k=0}^{[n/3]}\b{2k}k^2\b{4k}{2k}\b{n+k}{4k}(-3)^{n-3k}
\\&=\sum_{k=0}^n\b{2k}k^2\b{4k}{2k}\b{n+3k}{4k}(-27)^{n-k}.\endaligned
\tag 4.1$$
\pro{Theorem 4.1} Let $p$ be a prime with $p>3$. Then
$b_{p-1}\e 81^{p-1}\mod {p^3}$.
\endpro
Proof.  For $1\le k<\f p3$,
$$\align &(2k)!^2\b{4k}{2k}\b{p-1+k}{4k}
\\&=(p+k-1)(p+k-2)\cdots (p+1)p(p-1)(p-2)\cdots (p-3k)
\\&\e p(k-1)!(1+pH_{k-1})(-1)^{3k}(3k)!(1-pH_{3k})
\\&=p\cdot k!(3k)!\f{(-1)^k}k(1-p(H_{3k}-H_{k-1}))\mod {p^3}.
\endalign$$
Thus, appealing to (4.1),
$$\aligned \f{b_{p-1}}{3^{p-1}}
&=\sum_{k=0}^{[p/3]}\b{2k}k^2\b{4k}{2k}\b{p-1+k}{4k}\f 1{(-27)^k}
\\&\e 1+\sum_{k=1}^{[p/3]}\f{(-1)^k\cdot k!(3k)!}{k\cdot
k!^4}p(1-p(H_{3k}-H_{k-1}))\f 1{(-27)^k}
\\&=1+p\sum_{k=1}^{[p/3]}\f{\b{2k}k\b{3k}k}{k\cdot
27^k}-p^2\sum_{k=1}^{[p/3]}\f{\b{2k}k\b{3k}k}{k\cdot 27^k}
(H_{3k}-H_{k-1})\mod {p^3}.
\endaligned\tag 4.2$$
On the other hand, from (4.1) we have
$$\align
\f{b_{p-1}}{27^{p-1}}&=\sum_{k=0}^{p-1}\b{2k}k^2\b{4k}{2k}\b{p-1+3k}{4k}
\f 1{(-27)^k}
\\&=1+\sum_{k=1}^{p-1}\f{p(p-1)\cdots(p-k)\cdot
(p+1)(p+2)\cdots(p+3k-1)}{k!^4\cdot (-27)^k}
\\&\e 1+2p\sum_{k=[p/3]+1}^{p-1}\f{(-1)(-2)\cdots (-k)\cdot
1\cdot 2\cdots(3k-1)}{k!^4(-27)^k}\\&\q+p\sum_{k=1}^{[p/3]}
\f{(-1)(-2)\cdots(-k)(1-pH_k)\cdot
(3k-1)!(1+pH_{3k-1})}{k!^4(-27)^k}
\\&\e 1+2p\sum_{k=[p/3]+1}^{p-1}\f{(3k)!}{3k\cdot k!^3\cdot 27^k}
+p\sum_{k=1}^{[p/3]}\f{(3k)!}{k!^3\cdot 3k\cdot
27^k}(1+p(H_{3k-1}-H_k))
\\&=1+\f{2p}3\sum_{k=1}^{p-1}\f{\b{2k}k\b{3k}k}{k\cdot
27^k}-\f p3\sum_{k=1}^{[p/3]}\f{\b{2k}k\b{3k}k}{k\cdot 27^k}
+\f{p^2}3\sum_{k=1}^{[p/3]}\f{\b{2k}k\b{3k}k}{k\cdot
27^k}(H_{3k-1}-H_k)\mod {p^3} .\endalign$$ Hence
$$\align \f{b_{p-1}}{3^{p-1}}+3\cdot\f{b_{p-1}}{27^{p-1}}
&\e 1+p\sum_{k=1}^{[p/3]}\f{\b{2k}k\b{3k}k}{k\cdot
27^k}-p^2\sum_{k=1}^{[p/3]}\f{\b{2k}k\b{3k}k}{k\cdot 27^k}
(H_{3k}-H_{k-1})
\\&\q+3+2p\sum_{k=1}^{p-1}\f{\b{2k}k\b{3k}k}{k\cdot
27^k}-p\sum_{k=1}^{[p/3]}\f{\b{2k}k\b{3k}k}{k\cdot 27^k}
+p^2\sum_{k=1}^{[p/3]}\f{\b{2k}k\b{3k}k}{k\cdot
27^k}(H_{3k-1}-H_k)\\&=4+2p\sum_{k=1}^{p-1}
\f{\b{2k}k\b{3k}k}{k\cdot 27^k}-p^2\sum_{k=1}^{[p/3]}
\f{\b{2k}k\b{3k}k}{k\cdot 27^k}(H_{3k}-H_{3k-1}+H_k-H_{k-1})
\\&=4+2p\sum_{k=1}^{p-1}
\f{\b{2k}k\b{3k}k}{k\cdot 27^k}-p^2\sum_{k=1}^{[p/3]}
\f{\b{2k}k\b{3k}k}{k\cdot 27^k}\Big(\f 1{3k}+\f 1k\Big) \mod {p^3}
\endalign$$ and so
$$\Big(\f{1}{3^{p-1}}+\f 3{27^{p-1}}\Big)b_{p-1}
\e 4+2p\sum_{k=1}^{p-1} \f{\b{2k}k\b{3k}k}{k\cdot 27^k}-\f
43p^2\sum_{k=1}^{p-1} \f{\b{2k}k\b{3k}k}{k^2\cdot 27^k}\mod {p^3}
.\tag 4.3$$ Set $(x)_k=x(x+1)\cdots(x+k-1)$. It is evident that
$(x)_k=\b{-x}k(-1)^kk!$ and so $\f{(\f 13)_k(\f
23)_k}{(1)_k^2}=\b{-\f 13}k\b{-\f 23}k=\f{\b{2k}k\b{3k}k}{27^k}.$
Thus, from [T1] we deduce that
$$\sum_{k=1}^{p-1}
\f{\b{2k}k\b{3k}k}{k\cdot 27^k}\e -q_p\big(\f 1{27}\big)+\f
p2q_p\big(\f 1{27}\big)^2\mod{p^2}, \q\sum_{k=1}^{p-1}
\f{\b{2k}k\b{3k}k}{k^2\cdot 27^k}\e -\f 12q_p\big(\f
1{27}\big)^2\mod p.$$ It is easy to see that
$$\align q_p\big(\f
1{27}\big)&=\f{\ls
1{27}^{p-1}-1}p=\f{1-(1+pq_p(3))^3}{p(1+pq_p(3))^3}
\\&=-\f{3q_p(3)+3pq_p(3)^2+p^2q_p(3)^3}{1+3pq_p(3)+3p^2q_p(3)^2+p^3q_p(3)^3}
\e -\f {3q_p(3)+3pq_p(3)^2}{1+3pq_p(3)}
\\&\e -3q_p(3)(1+pq_p(3))(1-3pq_p(3))\e -3q_p(3)+6pq_p(3)^2\mod
{p^2}.\endalign$$ Hence
$$\align &\sum_{k=1}^{p-1}
\f{\b{2k}k\b{3k}k}{k\cdot 27^k}\e 3q_p(3)-6pq_p(3)^2+\f 92pq_p(3)^2=
3q_p(3)-\f 32pq_p(3)^2\mod {p^2},\tag 4.4
\\&\sum_{k=1}^{p-1}
\f{\b{2k}k\b{3k}k}{k^2\cdot 27^k}\e -\f 92q_p(3)^2\mod p.\tag
4.5\endalign$$  By (4.3)-(4.5),
$$\f{3+9^{p-1}}{27^{p-1}}b_{p-1}\e
4+2p(3q_p(3)-\f 32pq_p(3)^2)-\f 43p^2\Big(-\f 92q_p(3)^2\Big)
=4+6pq_p(3)+3p^2q_p(3)^2\mod {p^3}.$$ Since
$$ \f 1{3+9^{p-1}}=\f 1{4+2pq_p(3)+p^2q_p(3)^2} =\f{2-pq_p(3)}{2^3-p^3q_p(3)^3}\e \f
{2-pq_p(3)}8\mod {p^3},$$ From the above we deduce that
$$b_{p-1}\e 27^{p-1}\f 18(2-pq_p(3))(4+6pq_p(3)+3p^2q_p(3)^2)\e
27^{p-1}(1+pq_p(3))=81^{p-1}\mod {p^3},$$ which completes the proof.
\par\q
\par{\bf Remark 4.1} In [S11, Conjecture 4.1] the author
conjectured that for any prime $p>3$, $b_{p-1}\e 81^{p-1}-\f
2{27}p^3B_{p-3}\mod{p^4}$.
\par\q
\par
In [Su4], Z.W. Sun conjectured that for any prime $p\not=2,11$,
 $$\sum_{n=0}^{p-1}\b{2n}n\f{A_n'}{4^n}
 \e\cases x^2-2p\mod{p^2}&\t{if $\sls p{11}=1$ and so
 $4p=x^2+11y^2$,}\\0\mod{p^2}&\t{if $\sls p{11}=-1$.}
 \endcases$$
Now we partially solve the conjecture.
 \pro{Theorem 4.2} Let $p$ be
a prime with $p\not=2,11$. Then
$$\sum_{n=0}^{p-1}\b{2n}n\f{A_n'}{4^n}
 \e\cases x^2\mod p&\t{if $\sls p{11}=1$ and so
 $4p=x^2+11y^2$,}\\0\mod p&\t{if $\sls p{11}=-1$.}
 \endcases$$
 \endpro
Prof. Using Vandermonde's identity we see that for any nonnegative
integer $m$,
$$\align &\sum_{n=0}^m\b mn(-1)^{m-n}\sum_{k=0}^n\b nk^2\b{n+k}k
\\&=\sum_{n=0}^m\sum_{k=0}^n\b
mk\b{m-k}{n-k}\b{2k}k\b{n+k}{2k}(-1)^{m-n}
\\&=\sum_{k=0}^m\b
mk\b{2k}k(-1)^{m-k}\sum_{n=k}^m\b{m-k}{n-k}\b{n+k}{2k}(-1)^{n-k}
\\&=\sum_{k=0}^m\b
mk\b{2k}k(-1)^{m-k}\sum_{r=0}^{m-k}\b{m-k}r\b{2k+r}{r}(-1)^r
\\&=\sum_{k=0}^m\b
mk\b{2k}k(-1)^{m-k}\sum_{r=0}^{m-k}\b{m-k}{m-k-r}\b{-2k-1}r
\\&=\sum_{k=0}^m\b
mk\b{2k}k(-1)^{m-k}\b{m-3k-1}{m-k}
\\&=\sum_{k=0}^m\b
mk\b{2k}k\b{2k}{m-k} =\sum_{k=0}^m\b
mk\b{2(m-k)}{m-k}\b{2m-2k}k.\endalign$$ Note that $p\mid \b{2k}k$
for $\f p2<k<p$ and $\b{\f{p-1}2}k\e \b{2k}k4^{-k}\mod p$ for $k<\f
p2$. Taking $m=\f{p-1}2$ in the above we deduce that
$$\align \sum_{n=0}^{p-1}\b{2n}n\f{A_n'}{4^n}
&\e
(-1)^{\f{p-1}2}\sum_{n=0}^{(p-1)/2}\b{\f{p-1}2}n(-1)^{\f{p-1}2-n}A_n'
\\&=(-1)^{\f{p-1}2}\sum_{k=0}^{(p-1)/2}\b{\f{p-1}2}k
\b{2(\f{p-1}2-k)}{\f{p-1}2-k}\b{p-1-2k}k
\\&\e
(-1)^{\f{p-1}2}\sum_{k=0}^{(p-1)/2}\b{\f{p-1}2}k\b{\f{p-1}2}{\f{p-1}2-k}(-4)^{\f{p-1}2-k}
\b{-2k-1}k
\\&\e \sum_{k=0}^{(p-1)/2}\b{\f{p-1}2}k^24^{\f{p-1}2-k}\b{3k}k
\\&\e \sum_{k=0}^{(p-1)/2}\f{\b{2k}k^2\b{3k}k}{64^k}
\e  \sum_{k=0}^{p-1}\f{\b{2k}k^2\b{3k}k}{64^k}\mod p.
\endalign$$
Now applying [S5, Theorem 4.4] we deduce the result.
\par\q
\par By [M, (1.2)],
$$S_n=\sum_{k=0}^n\b{2k}k^2\b k{n-k}(-4)^{n-k}.$$
It is clear that
$$\align \sum_{n=0}^{p-1}c_n\f{S_n}{8^n}&
=\sum_{n=0}^{p-1}\f{c_n}{8^n}\sum_{k=0}^n\b{2k}k^2\b
k{n-k}(-4)^{n-k}
\\&=\sum_{k=0}^{p-1}\f{\b{2k}k^2}{(-4)^k}\sum_{n=k}^{p-1}\f{c_n}{(-2)^n}\b
k{n-k}
\\&=\sum_{k=0}^{p-1}\f{\b{2k}k^2}{(-4)^k}\sum_{r=0}^{p-1-k}\f{c_{k+r}}{(-2)^{k+r}}
\b kr
\\&=\sum_{k=0}^{p-1}\f{\b{2k}k^2}{8^k}
\sum_{r=0}^{p-1-k}\b kr\f{c_{k+r}}{(-2)^r}.\endalign$$ Thus,
$$\sum_{n=0}^{p-1}c_n\f{S_n}{8^n}
=\sum_{k=0}^{(p-1)/2}\f{\b{2k}k^2}{8^k} \sum_{r=0}^k\b
kr\f{c_{k+r}}{(-2)^r}+\sum_{k=(p+1)/2}^{p-1}
\f{\b{2k}k^2}{8^k}\sum_{r=0}^{p-1-k}\b kr\f{c_{k+r}}{(-2)^r}.\tag
4.6$$

 \pro{Theorem
4.3} Let $p$ be an odd prime. Then
$$\sum_{n=1}^{p-1}\f{nS_n}{8^n}\e 0\mod {p^2}.$$
\endpro
Proof. Observe that $p\mid \b{2k}k$ for $\f p2<k<p$. Taking $c_n=n$
in (4.6) gives
$$\align \sum_{n=1}^{p-1}\f{nS_n}{8^n}
&\e \sum_{k=0}^{(p-1)/2}\f{\b{2k}k^2}{8^k} \sum_{r=0}^k\b
kr\f{k+r}{(-2)^r}
\\&=\sum_{k=0}^{(p-1)/2}\f{\b{2k}k^2}{8^k}\Big(
k\sum_{r=0}^k\b kr\Big(-\f 12\Big)^r-\f k2\sum_{r=1}^k
\b{k-1}{r-1}\Big(-\f 12\Big)^{r-1}\Big)
\\&=\sum_{k=0}^{(p-1)/2}\f{\b{2k}k^2}{8^k}\Big(k\cdot\f 1{2^k}
-\f k2\cdot\f 1{2^{k-1}}\Big)=0\mod {p^2},\endalign$$ which proves
the theorem.
\par\q
\pro{Conjecture 4.1} Let $p$ be an odd prime. Then
$$\sum_{n=1}^{p-1}\f{nS_n}{8^n}\e\cases 0\mod {p^3}&\t{if $p\e 1\mod
4$,}
\\ 2p^2\mod {p^3}&\t{if $p\e 3\mod 4$.}
\endcases$$
\endpro
 \pro{Theorem 4.4} Let $p>3$ be a prime. Then
$$\sum_{n=0}^{p-1}\f{Q_n}{(-8)^n}\e 1\mod {p^2}\qtq{and}
\sum_{n=0}^{p-1}\f{Q_n}{(-9)^n}\e \Ls p3\mod {p^2}.$$
\endpro
Proof. Since $Q_n=\sum_{k=0}^n\b nk(-8)^{n-k}f_k$, using (2.1) we
see that
$$\align \sum_{n=0}^{p-1}\f{Q_n}{(-8)^n}&
=\sum_{n=0}^{p-1}\sum_{k=0}^n\b nk\f{f_k}{(-8)^k} =\sum_{k=0}^{p-1}
\f{f_k}{(-8)^k}\sum_{n=k}^{p-1}\b nk
\\&=\sum_{k=0}^{p-1}
\f{f_k}{(-8)^k}\b p{k+1}=\f{f_{p-1}}{(-8)^{p-1}}+\sum_{k=0}^{p-2}
\f{f_k}{(-8)^k}\cdot \f p{k+1}\b{p-1}k
\\&\e \f{f_{p-1}}{8^{p-1}}+p\sum_{k=0}^{p-2}\f{f_k}{(k+1)8^k}
=\f{f_{p-1}}{8^{p-1}}+p\sum_{k=1}^{p-1}\f{f_{p-1-k}}{(p-k)8^{p-1-k}}
\\&\e \f{f_{p-1}}{8^{p-1}}-p\sum_{k=1}^{p-1}\f{8^kf_{p-1-k}}k
 \mod{p^2}.\endalign$$ By [JV], $f_k\e (-8)^kf_{p-1-k}\mod p$.
Thus,
$$\sum_{n=0}^{p-1}\f{Q_n}{(-8)^n}\e \f{f_{p-1}}{8^{p-1}}-p\sum_{k=1}^{p-1}
(-1)^k\f{f_k}k\mod {p^2}.$$ By [Su3, Theorem 1.1 and Lemma 2.5],
$$f_{p-1}\e 1+3p\qp 2\mod {p^2}\qtq{and}\sum_{k=1}^{p-1}
(-1)^k\f{f_k}k\e 0\mod {p^2}.$$ Thus,
$$\sum_{n=0}^{p-1}\f{Q_n}{(-8)^n}\e \f{f_{p-1}}{8^{p-1}}\e \f
{1+3p\qp 2}{(1+p\qp 2)^3}\e 1\mod {p^2}.$$
\par Note that
$$\align &\sum_{k=0}^n\b nk(-9)^{n-k}a_k
\\&=\sum_{k=0}^n\b nk(-1)^{n-k}a_k\sum_{r=0}^{n-k}\b{n-k}r8^r
=\sum_{r=0}^n\b nr(-8)^r\sum_{k=0}^{n-r}\b{n-r}k(-1)^{n-r-k}a_k
\\&=\sum_{r=0}^n\b nr(-8)^rf_{n-r}=Q_n.\endalign$$ Using (2.1) we
see that
$$\align \sumnp\f{Q_n}{(-9)^n}&=\sumnp\sum_{k=0}^n\b
nk\f{a_k}{(-9)^k}=\sumkp\f{a_k}{(-9)^k}\sum_{n=k}^{p-1}\b nk
\\&=\sumkp\f{a_k}{(-9)^k}\b p{k+1}=\sumkp\f{a_k}{(-9)^k}
\cdot\f p{k+1}\b{p-1}k
\\&\e \f{a_{p-1}}{(-9)^{p-1}}+\sum_{k=0}^{p-2}\f{a_k}{9^k}
\cdot\f
p{k+1}=\f{a_{p-1}}{9^{p-1}}+\sum_{r=1}^{p-1}\f{a_{p-1-r}}{9^{p-1-r}}\cdot\f
p{p-r}
\\&\e \f{a_{p-1}}{9^{p-1}}-p\sum_{r=1}^{p-1}\f{9^ra_{p-1-r}}r\mod
{p^2}.\endalign$$ By [JV] or [S11, Theorem 3.1], $a_r\e \sls
p39^ra_{p-1-r}\mod p$. By [Su8, Lemma 3.2], $a_{p-1}\e \sls
p3(1+2p\qp 3)\e \sls p39^{p-1}\mod {p^2}$. Taking $x=1$ in [Su8,
(3.6)] gives $\sum_{r=1}^{p-1}\f {a_r}r\e 0\mod p$. Now, from the
above we deduce that
$$\align \sumnp\f{Q_n}{(-9)^n}&\e
\f{a_{p-1}}{9^{p-1}}-p\sum_{r=1}^{p-1}\f{9^ra_{p-1-r}}r \e \Ls
p3-p\Ls p3\sum_{r=1}^{p-1}\f {a_r}r\e \Ls p3\mod {p^2}.\endalign$$
This completes the proof.

\end{document}